\documentclass[sn-mathphys]{sn-jnl}
\jyear{2021}%

\theoremstyle{thmstyleone}%

\theoremstyle{thmstyletwo}%

\theoremstyle{thmstylethree}%

\usepackage{amsmath, amssymb}
\usepackage{latexsym}
\usepackage{amssymb}
\usepackage{amsmath}
\usepackage{txfonts}

\begin{document}

\title[Measure solutions to generalized Chaplygin gas equations]{Measure solutions to  piston problem for compressible fluid flow of  generalized chaplygin gas}

\author*[1]{\fnm{Meixiang} \sur{Huang}}\email{1228561480@qq.com}
\author[1]{\fnm{Yuanjin} \sur{Wang}}\email{454814482@qq.com}
\author[2]{\fnm{Zhiqiang} \sur{Shao}}\email{zqshao@fzu.edu.cn}

\affil*[1]{
\orgdiv{School of Mathematics and Statistics}, \orgname{Minnan Normal University}, \orgaddress{ \city{Zhangzhou}, \postcode{363000}, \state{Fujian}, \country{China}}}
\affil[2]{\orgdiv{College of Mathematics and Computer Science}, \orgname{ Fuzhou University}, \orgaddress{ \city{Fuzhou}, \postcode{350108}, \state{Fujian}, \country{China}}}

\abstract{We study the piston problem of the compressible fluid flow with the generalized Chaplygin gas. Depending on the inferential critical value of Mach number, we prove that, there exists an integral weak solution for the proceeding piston problem, consisting of a shock separating constant states ahead of the piston if Mach numbers less than this critical value, while a singular measure solution, with density containing a Dirac measure supported on the piston, shall be proposed to solve the proceeding piston problem if Mach numbers greater than or equal to the critical value. For the receding piston problem, rarefaction wave solution always exists when the piston recedes from the gas with any constant speed. Moreover, the  occurrence of vacuum state and  the convergence of solutions, as well as degeneration of equations are analyzed in the receding case as Mach number tends to infinity. }



\maketitle

\section{Introduction}\label{sec1}
We consider the Euler equations of one-dimensional compressible fluid flow in the following form:

\begin{equation}\label{key}
\left\{\begin{array}{ll} \rho_{t}+(\rho u)_{x}=0,\\
u_{t}+\left(\frac{u^{2}}{2}+\int_{}^{\rho}{\frac{P'(s)}{s}ds} \right) _{x}=0,
\end{array}\right . \tag{1.1}
\end{equation} 
where  $t\geq 0$ is a time variable, and  $x \in R$ is a space variable; $\rho $,  $u$ and $P$ represent the density of mass, the velocity and the pressure of the fluid, respectively. System (1.1) was first given by Earnshaw [1], and is regarded as the one-dimensional compressible fluid flow [2]. In addition, System (1.1) has other different physical backgrounds. For instance, it is a scaling limit system of Newtonian dynamics with long-range interaction for a continuous distribution of mass in $R$ [3] and also a hydrodynamic limit for the Vlasov equation [4].

There have been many studies on the Euler system of one-dimensional compressible fluid flow with different equation of state. Diperna [5] established the existence of global weak solution to the Cauchy problem for the case of $1<\gamma<3$ by using Glimm's scheme under the equation of state $P(\rho) = A \rho^{\gamma}$, where $A>0$ and $\gamma>1$. Based on this result, Li [6] solved the same problems in the case $-1<\gamma<1$. Lu [7] and  Cheng [8] proved analytically an existence theorem for the global entropy solutions in the case of $\gamma>3$ and $1<\gamma<3$ respectively, by means of the theory of compensated compactness coupled with some basic ideas of kinetic formulation. For the negative pressure laws,  Cheng et al. [9] studied the Riemann problem for the 1D compressible fluid flow of a Chaplygin gas and solved constructively the existence and uniqueness of delta shock solutions, and the model of generalized Chaplygin gas was also investigated by Pang et al [10]. 

In this work, we are concerned with  generalized Chaplygin gas,  the equation of state is as follows:

\begin{equation}\label{key}
P = - s\rho^{-\gamma},   ~~~~~~ 0 < \gamma \leq 1. \tag{1.2}
\end{equation}
 A distinctive feature of the generalized Chaplygin gas is that it has a negative pressure with the positive sound speed. The (1.2) with   $\gamma = 1$ is also called  the pure Chaplygin gas, which was introduced by Chaplygin [11], Tsien [12] and von Karman [13] as a mathematical approximation for calculating the lifting force on a wing of an airplane in aerodynamics. The generalized Chaplygin gas is usually used to describe the current accelerated expansion of the universe and the evolution of the perturbations of energy density [14,15]. 
 
 The piston problem, serving as an important typical physical model in mathematical fluid dynamics, has been studied extensively in the past decades [16-19].  Qu and Yuan [20] considered the piston problem for compressible Euler equations
 \begin{equation}\label{key}
\left\{\begin{array}{ll} \rho_{t}+(\rho u)_{x}=0,\\
(\rho u)_{t}+ (\rho u^{2}+P(\rho))_{x}=0
\end{array}\right . \tag{1.3}
  \end{equation}
with the Chaplygin gas.   They defined a Radon measure solution of piston problem within a convenient space of distributions that contains discontinuous functions and Dirac measures. Roughly speaking, the Radon measure solution is a solution such that at least one of the state variables has a Dirac delta function  with the boundary of the domain as its support.  Moreover, Gao, Qu and Yuan [21] established the equivalence of free piston and delta shock, for the one-dimensional pressureless compressible Euler equations, which helps to understand the physics of delta shocks, but also provides a way to reduce the fluid–solid interaction problem. Recently,  Fan et.al [22] investigated the piston problem for system (1.3) with the generalized Chaplygin gas (1.2), where the Radon measure solutions are  constructed to deal with the concentration of mass  on the piston. For the theory of piston problem and Radon measure solution,  readers can refer to [16, 20-23] for more details.

 To our knowledge, no literature  investigates the piston problem for model (1.1) with a negative pressure law (1.2). Motivated by the above discussions, we will study the piston problem of  system (1.1) and (1.2). First, let us describe the piston problem  as follows. Suppose there is a piston which can move leftward or rightward in an infinite long and thin tube extending along the horizontal x-axis. The tube is filled with gas, which is enclosed by the piston on the right hand side and is uniform and static initially. The goal of piston problem is to study what happens if the piston moves? 
 
For the convenience of treating this question, we insert (1.2) into (1.1) to give the following compressible fluid flow with the generalized Chaplygin gas
 
 \begin{equation}\label{key}
 \left\{\begin{array}{ll} \rho_{t}+(\rho u)_{x}=0,\\
 u_{t}+\left(\frac{u^{2}}{2}-A\rho^{-\alpha} \right) _{x}=0,
 \end{array}\right . \tag{1.4}
 \end{equation} 
 where $A=\frac{s\gamma}{1+\gamma}$, $\alpha = \gamma+1$. (1.4) is totally equivalent to (1.1) and (1.2).
 
 Suppose initially the piston lies at $x=0$, and the gas fills the domain $\{x<0\}$. We assume the gas is static and uniform with given state
 \begin{equation}\label{key}
U_{0} = (\rho, u) \vert_{t=0} = (\rho_{0}, 0). \tag{1.5}
\end{equation}
The piston moves with a given constant speed $v_{0}$ and the trajectory of the piston is $x = v_{0}t ~~(t\geq 0)$. Then, the time-space domain we consider is 
 \begin{equation}\label{key}
\Omega_{t} = \left\lbrace (t, x) : x< v_{0}t, ~~~t>0\right\rbrace. \tag{1.6}
\end{equation}
On the   trajectory of the piston, we impose the  impermeable condition
\begin{equation}\label{key}
 u(t, x) = v_{0}   ~~~~~    on~~~  x = v_{0}t. \tag{1.7}
\end{equation}
The piston problem considered in this paper is to find a solution of (1.4) in the domain $\Omega_{t}$, satisfying (1.5) and (1.7).

 Recall that the local sonic speed in generalized Chaplygin gas is given by
 \begin{equation}\label{key}
c = \sqrt{P'(\rho)} = \sqrt{sr}\rho^{-\frac{r+1}{2}}.  \tag{1.8}
\end{equation}
 Then, the Mach number $M_{0}$ of the piston with respect to the gas is defined by
 \begin{equation}\label{key}
M_{0} = \frac{ \vert v_{0} \vert }{c_{0}},\tag{1.9}
\end{equation}
where $v_{0}$ is  the move speed of piston, $c_{0}$ is as in (1.8) for $\rho = \rho_{0}$. In this paper, we will show that if the  piston moves into the gas with $M_{0} \in \left( 0, \sqrt{2}(1+\gamma)^{-\frac{1}{2}}\right) $, or recedes from the gas with any constant speed,  there exists integral weak solutions, consisting of shocks or rarefaction waves respectively. However, there is no piecewise constant integral weak solution if the piston moves into the gas with $M_{0} \geq \sqrt{2}(1+\gamma)^{-\frac{1}{2}}$, and then a concept of Radon measure solution is proposed  to understand the compressible fluid flow equations when the unknowns are measures, which also demenstrates the necessity of considering general measure solutions in the study of piston problems of systems of hyperbolic conservation laws.

The outline of this paper is organized as follows. In Section 2, we  reformulate the piston problem in a convenient way by shifting the coordinates system to move with the piston with Galilean transformation, and present a definition of measure solutions. The main results of piston problem are Theorem 2.1 and Theorem 2.2 given at the end of this section. In Section 3, we give detailed proof to  Theorem 2.1. First, we construct an integral weak solution consisting of a shock wave when the piston rushes into the gas and Mach number less than $\sqrt{2}(1+\gamma)^{-\frac{1}{2}}$. Then,  we find a solution containing a weighted Dirac measure supported on the piston if $M_{0}\geq\sqrt{2}(1+\gamma)^{-\frac{1}{2}}$, which also justified the concept of measure solutions we proposed. In Section 4, we rigorously prove Theorem 2.2. We  construct a self-similar solution containing rarefaction waves as the piston recedes from the gas, and find the occurrence of vacuum state when Mach number tends to infinity.

\section{The piston problem of generalized Chaplygin gas and main results}\label{sec2}

\subsection{Simplified piston problem }
In this part, we transform the piston problem of generalized Chaplygin gas  equations into simplified version by Galilean transformation. Specifically, 
 the following Galilean transformation is adopted  to shift the coordinates to move with the piston and reformulate the piston problem:

\begin{equation}\label{key}
\left\{\begin{array}{ll} 
t' = t,\\
x' = x - v_{0}t,\\
\rho'(t', x') = \rho(t', x'+v_{0}t'),\\
u'(t', x') = u(t', x'+v_{0}t')-v_{0},\\
P'(t', x') = P(t', x'+v_{0}t').
\end{array}\right. \tag{2.1} 
\end{equation}
\noindent
It is easy to validate that the equations  (1.1) are invariant under (2.1), and the domain $\Omega_{t}$ in (1.6) and  trajectory of the piston are reduced to $\Omega' = \left\lbrace (t', x') : x'< 0, t'>0\right\rbrace  $ and $x' = 0$. For convenience of statement, 
we drop all the primes "$'$" without confusion. So,  the domain  $\Omega'$ can be written by 

\begin{equation}\label{key}
\Omega = \left\lbrace (t, x) : x< 0, t>0\right\rbrace.  \tag{2.2}
\end{equation}
 Similarly, the initial condition becomes 
\begin{equation}\label{key}
(\rho, u)\vert_{t=0} = (\rho_{0}, -v_{0}),  \tag{2.3}
\end{equation}
 and the boundary condition becomes 
\begin{equation}\label{key}
 u(t, x) = 0   ~~~~~    on~~~  x = 0. \tag{2.4}
\end{equation}
\noindent
To better understand the essence of piston problem, we carry out the non-dimensional linear transformations of independent and dependent variables, which corresponds to some similarity laws in physics [20]:

\begin{equation}\label{key}
\tilde{t} = \frac{t}{T},~~~~
\tilde{x} = \frac{x}{L},~~~~
\tilde{\rho} = \frac{\rho}{\rho_{0}},~~~~
\tilde{u} = \frac{u}{\left \vert v_{0}\right \vert},~~~~
\tilde{P} = \frac{P}{\rho_{0}v_{0}^{2}}, \tag{2.5}
\end{equation}
where $T$ and   $L > 0$ are constants with $\frac{L}{T} = v_{0} $. Inserting (2.5) into (1.1) and calculating directly, we can obtain that

\begin{equation}\label{key}
\left\{\begin{array}{ll} \tilde{\rho}_{\tilde{t}} + (\tilde{\rho} \tilde{u})_{\tilde{x}} = 0,\\
\tilde{u}_{\tilde{t}}+\left(\frac{\tilde{u}^{2}}{2}+\int_{}^{\tilde{\rho}}{\frac{\tilde{P'}(s)}{s}ds} \right) _{\tilde{x}}=0,
\end{array}\right. \tag{2.6} 
\end{equation}
which, together with (2.5), implies that
 $\tilde{\rho}$, $\tilde{u}$, $\tilde{P}$ still satisfies (1.1) and hence (1.1) is  invariant under (2.5),  that is to say, the $\rho$, $u$ and $P$ in (1.1) are equivalent to $\tilde{\rho}$, $\tilde{u}$, $\tilde{P}$ in (2.6). So, from (2.5), in the following we shall take initial data as
\begin{equation}\label{key}
\rho_{0} = 1,~~ v_{0} = \pm 1. \tag{2.7}
\end{equation}
On the other hand, by using (1.8), substituting $c_{0} =  \sqrt{sr}\rho_{0}^{-\frac{\gamma+1}{2}}$ into $M_{0} = \frac{ \vert v_{0} \vert }{c_{0}}$, we  obtain
\begin{equation}\label{key}
s = \frac{\rho_{0}^{\gamma+1} v_{0}^{2}}{\gamma M_{0}^{2}}
.\tag{2.8}
\end{equation}
Then, (1.2) becomes
\begin{equation}\label{key}
P(\rho) = -\frac{\rho_{0}^{\gamma+1} v_{0}^{2}}{\gamma  M_{0}^{2}\rho^{\gamma}}. \tag{2.9}
\end{equation}
Combining (2.5) with (2.9) and replacing $\frac{\rho_{0}}{\rho}$ by $\frac{1}{\tilde{\rho}}$, we obtain 
\begin{equation}\label{key}
\tilde{P}(\tilde{\rho}) = -\frac{1}{  \gamma M_{0}^{2}  \tilde{\rho}^{\gamma} }. \tag{2.10}
\end{equation}
Since $\tilde{P}$, $\tilde{\rho}$  in (2.10) are equivalent to the corresponding $P$, $\rho$ in (2.9),  for abbreviation, we write

 \begin{equation}\label{key}
 P(\rho)= -\frac{1}{  \gamma M_{0}^{2}  \rho^{\gamma} }. \tag{2.11}
 \end{equation}
By (2.11) and taking $s= \frac{1}{\gamma M_{0}^2}$, we have $P(\rho) = -s \rho^{-\gamma}$, which also infers  that $\rho_{0} = 1$, $ v_{0} = \pm 1$. 
Therefore, we define  initial data by

\begin{equation}\label{key}
\rho_{0} = 1,~~ v_{0} = \pm 1,~~ P_{0} = -\frac{1}{\gamma M_{0}^{2}}. \tag{2.12}
\end{equation}
In conclusion, the piston problem can now be reformulated as to define and seek solutions of (1.4), (2.4) and (2.12) in the domain given by (2.2).

\subsection{Definition of measure solutions of piston problem}\label{sec3}
The above formulation of piston problem only makes sense for classical  elementary waves solutions such as shock waves or rarefaction waves. Since it turns out that the unknowns might be measures singular to Lebesgue measure, it is necessary to rewrite the piston problem and give a rigorous definition of Radon measure solution of compressible fluid flow for the mathematical analysis. 

Recall that a Radon measure $m$ on the upper plane $[0, \infty) \times R$ could act on the compactly supported continuous functions
\begin{equation}\label{key}
<m, \phi>  =  \int_{0}^{\infty}\int_{R}\phi(t,x)m(dxdt),  \tag{2.13}
\end{equation}
where the test function $\phi \in C_{0}([0, \infty)\times R)$. The standard  Lebesgue measure of $ R^{2}$ is denoted by $\mathcal{L}^{2}$. The fact that  a measure $\mu$  is absolutely continuous with respect to a nonnegative  measure $\nu$ is denoted by $\mu \ll \nu$. The Dirac measure supported on a curve, which is singular to $\mathcal{L}^{2}$ is defined as below [20].

\textbf{Definition 2.1.}~~Let $L$ be a Lipschitz curve given by $x=x(t)$ for  $t$ $  \in$ $ [0,T)$, and $w_L(t) \in L_{loc}^{1}(0,T)$. The Dirac measure $w_{L}\delta_{L}$ supported on $L \subset R^{2} $ with weight $w_{L}$ is defined by
\begin{equation}\label{key}
<w_{L}\delta_{L}, \phi> = \int_{0}^{T}\phi(t, x(t)) w_L(t) \sqrt{x'(t)^2+1}dt,~~~~~~\forall \phi \in C_{0}(R^2).  \tag{2.14}
\end{equation}

Now, we could formulate the piston problem rigorously by introducing the following definition of measure solutions.

\textbf{Definition 2.2.}~~ For fixed $0<M_{0}\leq \infty$, let $ \varrho,  m, n,  \tau, \mathcal{P} $ be Radon measures on $\overline{\Omega}$ (the closure of $\Omega$), and $w_{p}$ a locally integrable nonnegative function on $[0, \infty)$. We call $(\varrho, u, w_{p})$  a measure solution to the piston problem (1.4) (2.4) and (2.12), provided that

\noindent i) $m\ll \varrho$, $\tau \ll n$, and they have the same Radon-Nikodym derivative $u$; namely
\begin{equation}\label{key}
u \triangleq \dfrac{m(dxdt)}{\varrho(dxdt)} = \dfrac{\tau(dxdt)}{n(dxdt)};
\tag{2.15}
\end{equation}
ii) For any $\phi \in C_{0}^{1}(R^{2})$, there holds
\begin{equation}\label{key}
<\varrho, \partial t \phi> + <m, \partial x \phi> + \int_{-\infty}^{0}\rho_{0}\phi(0, x)dx = 0, \tag{2.16}
\end{equation}
\begin{equation}\label{key}
\begin{split}
&<n, \partial t \phi> + <\tau, \partial x \phi> +<\mathcal{P}, \partial x \phi> - <w_{p}\delta_{\{x=0, ~t \geq 0\} }, \phi>
+ \int_{-\infty}^{0}u_{0}\phi(0, x)dx = 0, 
\end{split}\tag{2.17}
\end{equation}
iii) If $\varrho \ll \mathcal{L}^2$ with derivative $\rho(t, x)$ in a neighborhood of $(t, x) \in [0, \infty) \times (-\infty, 0]$, and $\mathcal{P} \ll \mathcal{L}^2$ with derivative $P(t, x)$ there, then $\mathcal{L}^2$-a.e. there holds 
\begin{equation}\label{key}
P = -\frac{1}{ \rho^{\gamma}}\frac{1}{\gamma M_{0}^{2}}, \tag{2.18}
\end{equation}
and in addition, the classical entropy condition holds for discontinuities of functions $\rho, u$ near $(t, x)$.

\textbf{Remark 2.1.} Physically, the weight $w_{p}$ is the inertial force caused by the flows hitting the piston. It is alway positive when the piston moves to the gas with Mach number exceeding $\sqrt{2}(1+\gamma)^{-\frac{1}{2}}$, which will be rigorously confirmed in Section 3. Hence, compared with the standard pressureless Euler equations,  the high Mach number limit is not simply the vanishing pressure limit, since there is an extra term $w_{p}\delta_{\{x=0, t\geq 0\}}$ in the limiting Euler equations (see the fourth term in (2.17)). The requirement that $w_{p}$ is nonnegative shall be considered as a kind of stability condition for the measure solutions.

\textbf{Remark 2.2.} The above definition is consistent with the familiar  integrable weak solution just by taking
\begin{equation}
\varrho = \rho I_{\Omega}\mathcal{L}^2, ~~m = \rho u I_{\Omega}\mathcal{L}^2, ~~n =u I_{\Omega}\mathcal{L}^2,  ~~\tau = \frac{u^{2}}{2} I_{\Omega}\mathcal{L}^2, \mathcal{P} = \int_{}^{\rho}\frac{P'(s)}{s}ds I_{\Omega}\mathcal{L}^2. \notag
\end{equation}
The main results of this paper are the following two theorems.

 \textbf{Theorem 2.1.} ~~ For the piston moving into the gas ($v_{0} = -1$), there has two situations. One is that a shock wave solution of problem (1.1), (2.4) and (2.12) exists only for Mach number satisfying $M_{0} \in \left( 0, \sqrt{2}(1+\gamma)^{-\frac{1}{2}}\right) $. The other is that  a Radon measure solution we defined is justified as a solution to the problem (1.1), (2.4) and (2.12) in domain  $\Omega$  if Mach number $M_{0} \geq \sqrt{2}(1+\gamma)^{-\frac{1}{2}}$.

 \textbf{Theorem 2.2.} \label{{Theorem 4.2}}~~ For  the piston recedes from the gas ($v_{0} = 1$), the problem (1.1), (2.4) and (2.12) always has a rarefaction wave solution in the domain $\Omega$ for $0<M_{0}<\infty$. If $M_{0} \rightarrow 0$, the vacuum state occurs in the front of the piston and the rarefaction wave degenerates the discontinuity line $x = - t$, the high Mach number limiting equations are consistent with the pressureless Euler equations.

In what follows, we complete the detailed proofs for these two theorems in Section 3 and 4.
\section{Proof of Theorem 2.1} 
\subsection{Shock wave solution for $M_{0} \in \left( 0, \sqrt{2}(1+\gamma)^{-\frac{1}{2}}\right) $} 
~~~~Since shock wave appears ahead of the piston when it pushes into the gas, we first give the definition of the integral weak solution for problem (1.1), (2.4) and (2.12).

 \textbf{Definition 3.1.}~~ We say $\left(\rho, u \right)\in L^{\infty}(\left[ 0,\infty\right) \times\left( -\infty,0\right] ) $ is an integral solution to  the problem (1.1), (2.4) and (2.12), if for any $\phi \in C^{1}_{0}(R^{2})$, there hold

\begin{equation}\label{key}
\left\{\begin{array}{ll} \int_{\Omega}(\rho \partial t \phi+\rho u \partial x \phi)dxdt + \int_{-\infty}^{0}\rho_{0}(x)\phi(0,x)dx = 0,\\
\\
\int_{\Omega}\left[   u \partial t \phi +\left(  \frac {u^{2}}{2}+\int_{}^{\rho}{\frac{P'(s)}{s}ds}\right)  \partial x \phi \right]  dxdt \\-
\int_{0}^{\infty}\int_{}^{\rho}{\frac{P'(s)}{s}ds}\big\vert_{x=0} \phi(t,0)dt 
+ \int_{-\infty}^{0}u_{0}(x) \phi(0,x)dx = 0.
\end{array}\right.
\tag{3.1}
\end{equation}

Notice that the problem (1.1), (2.4) and (2.12) is a Riemann problem with boundary conditions taken into account, and the equations (1.1), the initial data (2.12) as well as the boundary condition (2.4) are invariant under the scaling 
$$
(x, t)\rightarrow (\mu x, \mu t)   ~~~~~for ~~~\mu  \neq 0.
$$
Thus, we can construct a piecewise constant
self-similar solution $U(x, t) = V(\frac{x}{t}) $ to connect two states $(\rho_{0}, -v_{0})$ and $(\rho_{1}, 0)$, which is in the form

\begin{equation}\label{key}
	U(x, t) = V\left( \frac{x}{t}\right) =
	\left\{\begin{array}{ll} 
		V_{0} = (1,1),~~~~~ -\infty \leq \frac{x}{t} < \sigma,\\
		V_{1} = (\rho_{1}, 0),~~~~~  \sigma < \frac{x}{t} \leq 0, \tag{3.2}
	\end{array}\right .
\end{equation}
where $V_{0}$ and $V_{1}$ are subject to the Rankine-Hugoniot condition

\begin{equation}\label{key}
\left\{\begin{array}{ll} 
\sigma(\rho_{1} - \rho_{0}) = \rho_{1}u_{1} - \rho_{0}u_{0},~~~~\\
\sigma(u_{1} -u_{0}) = \frac{u_{1}^{2}}{2} -A\rho_{1}^{-\alpha} - \frac{u_{0}^{2}}{2} +A\rho_{0}^{-\alpha},
\end{array}\right. \tag{3.3}
\end{equation}
where $A=\frac{sr}{1+r}$, ~~$\alpha = \gamma+1$.
In view of $\rho_{0} = 1$, $u_{0} = 1$, $u_{1} = 0$, it follows from the first equation of (3.3) that

\begin{equation}\label{key}
\sigma = -\frac{1}{\rho_{1}-1}. \tag{3.4}
\end{equation}
\noindent
Note that $\sigma <0$ requires that $\rho_{1} >1$. Substituting $\sigma$ into the second equation of (3.3) gives

\begin{equation}\label{key}
\frac{1}{\rho_{1}-1}
= - \frac{s\gamma \rho_{1}^{-\alpha}}{1+\gamma}+\frac{s\gamma }{1+\gamma}-\frac{1}{2}. \tag{3.5}
\end{equation}
where $s = \frac{1}{\gamma M_{0}^{2}}$.  Simplifying and rearranging (3.5) gives that
\begin{equation}\label{key}
M_{0}^2 = \frac{2(1-\rho_{1}^{-1})(1-\rho_{1}^{-\alpha})}{(1+\gamma)(1+\rho_{1}^{-1})} = \frac{2(1-\rho_{1}^{-1})(1-\rho_{1}^{-\gamma-1})}{(1+\gamma)(1+\rho_{1}^{-1})}
. \tag{3.6}
\end{equation}
For $0<M_{0}<\sqrt{2}(1+\gamma)^{-\frac{1}{2}}$, we have $0< \frac{(1-\rho_{1}^{-1})(1-\rho_{1}^{-\gamma-1})}{(1+\rho_{1}^{-1})}<1$. 
Define a continuous function  $f(\rho)= \frac{(1-\rho^{-1})(1-\rho^{-\gamma-1})}{(1+\rho^{-1})}$ with respect to $\rho$. Then, we can deduce $f(1) = 0$, $f(\infty) =1 $. According to the intermediate value theorem of continuous function,  there exists a $\rho_{1} > 1$, such that $0 < f(\rho_{1}) < 1$, which shows the existence of a shock wave solution. Furthermore, a direct computation shows that
\begin{equation}\label{key}
\begin{split}
f'(\rho) &=  \frac{\left[ 2\left( 1-\rho^{-\gamma-1}\right)+2(\rho-1)(\gamma+1) \rho^{-\gamma-2} 
	\right](\rho+1)-2(\rho-1)\left( 1-\rho^{-\gamma-1}\right) }
{(1+\gamma)(\rho+1)^2}
\\&=\frac{4\left( 1-\rho^{-\gamma-1}\right)+2(\rho-1)(\rho+1)(\gamma+1) \rho^{-\gamma-2}  }
{(1+\gamma)(\rho+1)^2}.
\end{split}\tag{3.7}
\end{equation}
From (3.7), it is easy to infer $f'(\rho)>0$ for $\rho>1$, i.e., $f(\rho)$ is always a monotone and increasing function for $\rho>1$, 
which guarantees the uniqueness of the shock wave solution. Thus, we rigorously prove that there has a shock wave solution (3.2) connecting two states $(1,1)$ and $(\rho_{1}, 0)$ as the piston rushes into the gas with Mach number $0<M_{0}<\sqrt{2}(1+\gamma)^{-\frac{1}{2}}$.

\subsection{Singular measure solution for $M_{0}\geq\sqrt{2}(1+\gamma)^{-\frac{1}{2}}$}

From (3.6), if $M_{0}\geq\sqrt{2}(1+\gamma)^{-\frac{1}{2}}$, it follows that $f(\rho_{1})\geq 1$, which contradicts with $\rho_{1}>1$. In this case, the piston problem cannot be solved in the sense of integral weak solutions which are measurable functions with respect to the Lebesgue measure.   To solve this case, we construct a  special measure solution by supposing that
 \begin{equation}\label{key}
\varrho =   I_{\Omega}\mathcal{L}^2 + w_{\rho}(t)\delta_{\{x=0,~ t \geq 0\} }, 
~m =   I_{\Omega}\mathcal{L}^2 , ~n =  I_{\Omega}\mathcal{L}^2,   ~\tau =  \frac{1}{2} I_{\Omega}\mathcal{L}^2,~\mathcal{P} = -\frac{1}{(1+\gamma)M_{0}^2}I_{\Omega}\mathcal{L}^2. \tag{3.8}
\end{equation}
where $ I_{\Omega}$ is the characteristic function of $\Omega$. These expressions come from the physical phenomena of infinite-thin shock layer in hypersonic flows past bodies and the hypersonic similarity law [24], and observations made in [17, Remark 1].

According to the Definition 2.2 and integration-by-parts, we obtain
\begin{equation}\label{key}
\begin{split}
0 &= <\varrho, \partial t \phi> + <m, \partial x \phi> + 
\int_{-\infty}^{0} \rho_{0} \phi(0, x)dx \\&= \int_{\Omega} \partial t \phi dxdt + \int_{0}^{\infty} w_{\rho}(t) \partial t  \phi(t, 0)dt
+   \int_{0}^{\infty}  \int_{-\infty}^{0}  \partial x  \phi dxdt 
+ \int_{-\infty}^{0} \phi(0, x)dx 
\\&= \int_{-\infty}^{0} \phi(t, x)\vert_{t=0}^{\infty} dx
+ w_{\rho}(t) \phi(t, 0)\vert_{t=0}^{\infty} 
-  \int_{0}^{\infty} w'_{\rho}(t) \phi(t, 0)dt \\&+ 
  \int_{0}^{\infty} \phi(t, x)\vert_{x=-\infty}^{0} dt +  \int_{-\infty}^{0} \phi(0, x)dx 
\\&=  \int_{0}^{\infty} (1 - w'_{\rho}(t))\phi(t, 0)dt
-w_{\rho}(0)\phi(0, 0).\end{split} \tag{3.9}
\end{equation}
By the arbitrariness of $\phi$, it follows
\begin{equation}\label{key}
\left\{\begin{array}{ll} 
w'_{\rho}(t) = 1, ~~~t\geq0,\\
w_{\rho}(0) = 0.
\end{array}\right. \tag{3.10}
\end{equation}
Thus, we have 
\begin{equation}\label{key}
w_{\rho}(t) = t. \tag{3.11}
\end{equation}
From the momentum equation (2.17), similar calculations show that
\begin{equation}\label{key}
\begin{split}
0 &= <n, \partial t \phi> + <\tau, \partial x \phi> + 
<\mathcal{P}, \partial x \phi> - <w_{p}(t)\delta_{\{x=0,~ t \geq 0\} }, \phi>
+ \int_{-\infty}^{0} u_{0} \phi(0, x)dx 
\\&=  \int_{0}^{\infty}  \int_{-\infty}^{0}
\partial t \phi dxdt +
\frac{1}{2} \int_{0}^{\infty}  \int_{-\infty}^{0}
\partial x \phi dxdt - 
\frac{1}{(1+\gamma)M_{0}^2} \int_{0}^{\infty}  \int_{-\infty}^{0}
\partial x \phi dxdt \\&-
\int_{0}^{\infty} w_{p}(t)  \phi(t, 0)dt
+    \int_{-\infty}^{0} \phi(0, x)dx 
\\&= \int_{0}^{\infty} \left[ \frac{1}{2}-w_{p}(t) - \frac{1}{(1+\gamma)M_{0}^2}\right] \phi(t, 0)dt.
\end{split} \tag{3.12}
\end{equation}
Since the arbitrariness of $\phi$ and $M_{0}\geq \sqrt{2}(1+\gamma)^{-\frac{1}{2}}$, then
\begin{equation}\label{key}
w_{p}(t) = \frac{1}{2} -\frac{1}{(1+\gamma)M_{0}^2}  \geq 0. \tag{3.13}
\end{equation}
Therefore, 
 the shock wave and all gas between the shock wave and  the piston adhere to the piston and then form a concentration of mass like a Dirac measure  as $M_{0}\geq\sqrt{2}(1+\gamma)^{-\frac{1}{2}}$.  From Remark 2.1., we also see that the high Mach number limit is  different from the the vanishing pressure limit. As the Mach number $M_{0}$ increases to infinity, the shock-front approaches the surface of the piston. 
 
 This completes the proof of Theorem 2.1.

\section{Proof of Theorem 2.2}

~~~~This section is devoted to the proof of Theorem 2.2.  Let $U(x, t) = (\rho(x, t), u(x, t))$. The system (1.4) has two eigenvalues [10]

 \begin{equation}\label{key}
 \lambda_{1}(U) = u -\sqrt{A\alpha}\rho^{-\frac{\alpha}{2}}, ~~~~~ \lambda_{2}(U) = u +\sqrt{A\alpha}\rho^{-\frac{\alpha}{2}}, \tag{4.1}
 \end{equation}
 with two associated right eigenvectors
 
 \begin{equation}\label{key}
 \vec{r}_{1} = \left(  \sqrt{\rho}, -\sqrt{A\alpha \rho^{-(\alpha+1)}}\right)^{T},~~~~~
  \vec{r}_{2} =  \left(  \sqrt{\rho}, \sqrt{A\alpha \rho^{-(\alpha+1)}}\right)^{T}, \tag{4.2}
 \end{equation}
satisfying 
\begin{equation}\label{key}
\triangledown \lambda_{i}(U) \cdot\vec{r}_{i} = \mp \frac{2-\alpha}{2}
\sqrt{A\alpha \rho^{-(\alpha+1)}}
 \neq 0, ~~~ (i=1, 2), \tag{4.3}
\end{equation}
where $A= \frac{s\gamma}{1+\gamma}$, $\alpha=\gamma+1$, $s=\frac{1}{\gamma M_{0}^2}$. Thus, all of the characteristic fields are genuinely nonlinear and  the Riemann solutions of the system (1.4) contain the first family rarefaction wave $R_{1}(U)$ and the second family rarefaction wave $R_2(U)$.

Now, we can as well construct a solution of the form $U(t,x) = V\left( \frac{x}{t}\right) $. For any fixed $0<M_{0}<\infty$, we suppose the solution is composed of two constant states $V_{0} = (1, -1)$, $V_{1} = (\rho_{1}, 0)$, and a rarefaction wave $V_{m}$ connecting them:

\begin{equation}\label{key}
	U(x, t) = V\left( \frac{x}{t}\right)  =
\left\{\begin{array}{ll} 
V_{0} ,~~~~~~~~~~~~~~~ -\infty \leq \frac{x}{t} < \lambda_{i}(V_{0}),\\

V_{m}\left(\frac{x}{t} \right)  ,~~~~~~~~~~ \lambda_{i}(V_{0}) \leq \frac{x}{t} < \lambda_{i}(V_{1}),\\

V_{1},~~~~~~~~~~~~~~~~~~  \lambda_{i}(V_{1}) < \frac{x}{t} \leq 0,
\end{array}\right. ~~~ (i=1, 2).  \tag{4.4}
\end{equation}

\noindent If the piston recedes from the gas, then $v_{0} >0$ and $\rho_{1} < \rho_{0}$. For the first family rarefaction wave $R_{1}(U)$, we can deduce the self-similar solution  as follow: 

\begin{equation}\label{key}
\left\{\begin{array}{ll} 
\eta = \frac{x}{t} = \lambda_{1}(U) = u - \sqrt{A\alpha}\rho^{-\frac{\alpha}{2}},\\
\\
u -\frac{2}{\alpha}\sqrt{A\alpha}\rho^{-\frac{\alpha}{2}} =
 u_{0} - \frac{2}{\alpha}\sqrt{A\alpha}\rho_{0}^{-\frac{\alpha}{2}} , ~~\rho_{1} < \rho < \rho_{0}, \\
\\
\lambda_{1}(V_{0}) \leq \lambda_{1}(U) \leq \lambda_{1}(V_{1}),
\end{array}\right.\tag{4.5}
\end{equation}
where $U(t, x) = V(\frac{x}{t})$. Replacing $\rho_{0}, u_{0}$ by $\rho_{0}=1, u_{0}=-1$ and the second equation of (4.5), we calculate that 

\begin{equation}\label{key}
u = -1+\frac{2}{\alpha}\sqrt{A\alpha}\left( \rho^{-\frac{\alpha}{2}}-1\right) . \tag{4.6}
\end{equation}
Then, substituting $A = \frac{s\gamma}{1+\gamma},~ \alpha=1+\gamma,~ s=\frac{1}{\gamma M_{0}^2}$ into (4.6) gives

\begin{equation}\label{key}
u = -1 + 2 (\gamma+1)^{-1}M_{0}^{-1}\left( \rho^{-\frac{\gamma+1}{2}}-1\right) .  \tag{4.7}
\end{equation}
Taking (4.7) into the first equation of (4.5), one obtains
\begin{equation}\label{key}
\eta = -1 - 2 (\gamma+1)^{-1}M_{0}^{-1}+M_{0}^{-1}(1-\gamma)(1+\gamma)^{-1} \rho^{-\frac{\gamma+1}{2}}.  \tag{4.8}
\end{equation}
where $0<\gamma<1$ and $0<M_{0}<\infty$. The density $\rho$ will be regard as a function of $\eta$. From (4.8), we have
\begin{equation}\label{key}
\rho(\eta) = \left[ \frac{(\eta+1)(1+\gamma)M_{0}+2}{1-\gamma}
\right]^{-\frac{2}{1+\gamma}} .  \tag{4.9}
\end{equation}
Combining (4.7) with (4.9), we obtain

\begin{equation}\label{key}
u(\eta) =-1+2(\eta+1)(1-\gamma)^{-1}+2M_{0}^{-1}(1-\gamma)^{-1}.  \tag{4.10}
\end{equation}
Considering $V_{0}=(1, -1)$, by the first equation of (4.5), we obtain $\eta_{-1} = \lambda_{1}(V_{0}) = 1-M_{0}^{-1}$. Moreover, for state $V_{1} = (\rho_{1}, 0)$, $u(\eta)$ satisfies the boundary condition $u(\eta_{0})=0$, which deduces $\eta_{0}=-\frac{\gamma+1}{2}-\frac{1}{M_{0}}$. Based on the above discussion and (4.9), we conclude that there has $V_{m}\left(\frac{x}{t} \right)=(\rho(\eta), u(\eta))$, satisfying $$\eta_{-1}\leq\eta\leq\eta_{0}$$ and $$0<\rho(\eta_{0})\leq\rho(\eta)\leq\rho(\eta_{-1}).$$ Consequently, (4.4) gives a rarefaction wave solution for system (1.1) with (2.4) and (2.12) in the domain $\Omega$.

If the Mach number $M_{0}\rightarrow \infty$, we obtain $\eta_{-1}\rightarrow -1$, $\eta_{0} \rightarrow -\frac{1+\gamma}{2}$, $P_{0} \rightarrow 0$, $P_{1} \rightarrow 0$. 

For arbitrarily given $\eta$ satisfying $\eta_{-1}<\eta\leq\eta_{0}$. If $M_{0} \rightarrow \infty$, we have $-1<\eta\leq-\frac{1+\gamma}{2}$. By this and (4.9), we obtain 
\begin{equation}\label{key}
\lim\limits_{M_{0}\rightarrow \infty}\rho(\eta)=\lim\limits_{M_{0}\rightarrow \infty} \left[ \frac {1}{(\eta+1)(1+\gamma)M_{0}+2}
\right]^{\frac{2}{1+\gamma}}\left(1-\gamma \right)^{\frac{2}{1+\gamma}}= 0.  \tag{4.11}
\end{equation}
Substituting $u(\eta_{0})=0$ into (4.7) yields $\rho(\eta_{0}) = \left[
\frac{M_{0}(\gamma+1)+2}{2}
 \right]^{-\frac{2}{\gamma+1}} $, it is easy to get
 
 \begin{equation}\label{key}
 \lim\limits_{M_{0}\rightarrow \infty}\rho(\eta_{0})=0.  \tag{4.12}
 \end{equation}
Based on above analysis and  (4.11) and (4.12), we conclude that the vacuum state occurs in the front of the piston and the
 right discontinuity line $x= - \frac{1+\gamma}{2}t$ vanishes if $M_{0}\rightarrow \infty$. Moreover, the rarefaction wave $\eta = \frac{x}{t} \in [\eta_{-1}, \eta_{0}]$
 degenerates discontinuity line $x=-t$, the limit equations for  high Mach number  are pressureless Euler equations. 

Analogously, for the second rarefaction wave $R_{2}(U(x,t))$, 
 we can deduce the self-similar solution  as follow: 

\begin{equation}\label{key}
\left\{\begin{array}{ll} 
\eta = \frac{x}{t} = \lambda_{2}(U) = u + \sqrt{A\alpha}\rho^{-\frac{\alpha}{2}},\\
\\
u + \frac{2}{\alpha}\sqrt{A\alpha}\rho^{-\frac{\alpha}{2}} =
u_{0} + \frac{2}{\alpha}\sqrt{A\alpha}\rho_{0}^{-\frac{\alpha}{2}} , ~~\rho_{1} < \rho < \rho_{0}, \\
\\
\lambda_{2}(V_{0}) \leq \lambda_{2}(U) \leq \lambda_{2}(V_{1}).
\end{array}\right.\tag{4.13}
\end{equation} 
Similar calculations lead to 
 \begin{equation}\label{key}
\rho(\eta) = \left[\left(\eta+1-2(1+\gamma)^{-1}M_{0}^{-1} \right)(1+\gamma) (\gamma-1)^{-1}M_{0}
 \right]^{-\frac{2}{1+\gamma}}, \tag{4.14}
\end{equation}
and
 \begin{equation}\label{key}
u(\eta) = -1+2(\eta+1)(1-\gamma)^{-1}+2(\gamma-1)^{-1}M_{0}^{-1}, \tag{4.15}
\end{equation}
where $-1+M_{0}^{-1}\leq \eta \leq -\frac{\gamma+1}{2}+M_{0}^{-1}$. Since $u_{1}=0$, we have
 \begin{equation}\label{key}
\rho_{1} = \left( 1-\frac{(1+\gamma)M_{0}}{2}\right)^{-\frac{2}{1+\gamma}} .  \tag{4.16}
\end{equation}
There are several different results. If $0<M_{0}<\frac{2}{1+\gamma}$, then $\rho_{1}>1$. Moreover, we have $\rho_{1}<0$ for $0<\gamma<1$ and $M_{0}>\frac{2}{1+\gamma}$. If $M_{0} \rightarrow \frac{2}{1+\gamma}$, we have
 \begin{equation}\label{key}
 \lim\limits_{M_{0}\rightarrow \frac{2}{1+\gamma}}
\left( 1-\frac{(1+\gamma)M_{0}}{2}\right)^{-\frac{2}{1+\gamma}}=\infty .  \tag{4.17}
\end{equation}
The above result contradicts with  the requirement of $0\leq \rho_{1}<\rho_{0}=1$. Therefore, the second rarefaction wave
$R_{2}(U(x,t))$ is not a physical solution as the piston recedes from the gas.

\textbf{Remark 4.1.} we remark that the limiting equations for compressible fluid flow of generalized chaplygin gas under the Mach number
$M_{0}\rightarrow \infty$ is the pressureless Euler equations. The rarefaction wave solutions of compressible fluid flow equations for piston problems at $x=-t$ degenerate into contact discontinuity.

\section{Acknowledgements}\label{sec6}

This research is supported by the Department of Education  of Fujian Province in China (No. JAT210254) and the  Minnan Normal University (No. KJ2021020).

 \section{ Statements and Declarations}
  \subsection{  Competing Interests}
The authors have no financial or non-financial interests related to this work to disclose.

\bibliography{sn-bibliography}

\end{document}